\documentclass[12pt]{article}
\usepackage{amscd, amssymb, latexsym, amsmath}

\newtheorem{pro}{Proposition}
\newtheorem{lem}{Lemma}

\newenvironment{proof}
{\noindent {\em Proof.}}
{\hfill $\Box$}

\numberwithin{thm}{section}
\numberwithin{cor}{section}
\numberwithin{pro}{section}
\numberwithin{lem}{section}
\numberwithin{dfn}{section}
\numberwithin{rem}{section}
\numberwithin{equation}{section}

\newcommand{\R}{\mathbb R}

\newcommand{\pk}{\partial_k}

\newcommand{\pj}{\partial_j}

\begin{document}
\title{Long-time Existence and Convergence of Graphic Mean Curvature Flow in
Arbitrary Codimension}
\author{Mu-Tao Wang}
\date{December 30, 2000, revised October 4, 2001}
\maketitle

\centerline{email:\,mtwang@math.columbia.edu}

\begin{abstract}
Let $f:\Sigma_1\mapsto \Sigma_2$ be a map between compact
Riemannian manifolds of constant curvature. This article considers
the evolution of the graph of $f$ in $\Sigma_1\times \Sigma_2$ by
the mean curvature flow. Under suitable conditions on the
curvature of $\Sigma_1$ and $\Sigma_2$ and the differential of the
initial map, we show that the flow exists smoothly for all time.
At each instant $t$, the flow remains the graph of a map $f_t$ and
$f_t$ converges to a constant map as $t$ approaches infinity. This
also provides a regularity estimate for Lipschtz initial data.
\end{abstract}

\section{Introduction}
The deformation of maps between Riemannian manifold has been
studied for a long time. The idea is to find a natural process to
deform a map to a "canonical" one. The harmonic heat flow is
probably the most famous example. It is the gradient flow of the
energy functional of maps. The classical work of Eells and Sampson
\cite{es} proves the flow converges to a harmonic map if the
target manifold is of non-positive curvature. However,
singularities do occur in the positive curvature case. Such
example exists even for maps between two-spheres. In \cite{ke},
\cite{lg}, the author proposes the study of a new deformation
process given by the mean curvature flow. The idea is to consider
the graph of the map as a submanifold in the product space and
evolve it in the direction of its mean curvature vector. This is
the gradient flow of the volume functional of graphs and a
stationary point is a "minimal map" first introduced by Schoen in
\cite{ sch}. This proposal proves to be quite successful in the
surface case. In fact, it provides a analytic proof of Smale's
classical theorem of the diffeomorphism group of
 spheres.
 This article considers the  arbitrary
  dimension and codimension case. We first prove that the graph condition is
preserved and the solution exists for all time.

\vskip 10pt \noindent {\bf Theorem A.} {\it Let $(\Sigma_1, g)$
and $(\Sigma_2, h)$ be Riemannian manifolds of constant curvature
$k_1$ and $k_2$ respectively
 and $f$ be
a smooth map from $\Sigma_1$ to $\Sigma_2 $. Suppose $k_1\geq
|k_2|$. If $\,\det(g_{ij}+(f^*h)_{ij})<2$, the mean curvature flow
of the graph of $f$
 remains a graph and exists for all time.} \vskip 10pt

The mean curvature flow for graphs appears to favor positively
curved domain manifold. The convergence theorem is the following.

\vskip 10pt \noindent {\bf Theorem B. } {\it  Let $(\Sigma_1, g)$
and $(\Sigma_2, h)$ be Riemannian manifolds of constant curvature
$k_1$ and $k_2$ respectively
 and $f$ be
a smooth map from $\Sigma_1$ to $\Sigma_2 $.
Suppose $k_1\geq |k_2|$ and $k_1+k_2
>0$. If $\,\det(g_{ij}+(f^*h)_{ij})<2$, then the mean curvature
flow of the graph of $f$ converges to the graph of a constant map
at infinity.} \vskip 10pt

The condition $\,\det(g_{ij}+(f^*h)_{ij})<2$ is actually a
geometric condition. When $\Sigma_1$ and $\Sigma_2$ are of the
same dimension, it is closely related but slightly stronger than
the condition that the Jacobian $J_1$ of the projection from the
graph to $\Sigma_1$ be strictly greater than the absolute value of
the Jacobian $J_2$ of the projection from the graph to $\Sigma_2$.
The geometric meaning is that we see more of the graph from
$\Sigma_1$ than from $\Sigma_2$.

An assumption of this type is clearly needed in general. In the
two dimensional case , the condition $J_1>|J_2|$ turns out to be
the optimal one in \cite{ke}. But in higher dimension, is is not
yet clear what would be the optimal condition for the global
existence of the flow.

 The stability of general gradient flow on
Riemannian manifolds with analytic metric was proved by Simon in
\cite{si} where the smallness of the second derivative is assumed
initially. The assumption in Theorem B is a condition on first
derivatives and the proof relies on a curvature estimate which
implies regularity for initial data with small Lipschitz norm.
Such estimate for codimension one graphic mean curvature flow was
proved by Ecker and Huisken \cite{eh}. A localized version in
codimension one case indeed gives an alternative proof of the
short time existence of the mean curvature flow.

Theorem A and B are proved in \S 4. The author would like to thank
Professor R. Schoen, Professor L. Simon and Professor S.-T. Yau
for their encouragement and advice.

\section{Preliminaries}

Let $f:\Sigma_1\mapsto \Sigma_2$ be a smooth
map between Riemannian manifolds.
The graph of $f$ is an embedded submanifold $\Sigma$
in $M=\Sigma_1\times \Sigma_2$.
We denote the embedding by $F:\Sigma_1
\mapsto M$, $F=id \times f$. There are
 isomorphisms $T\Sigma_1\mapsto T\Sigma$
by $X\mapsto X+df(X)$
and $T\Sigma_2\mapsto N\Sigma$
by $Y \mapsto Y-(df)^T(Y)$ where
$(df)^T: T\Sigma_2 \mapsto T\Sigma_1$
is the adjoint of $df$.

We assume the mean curvature flow of $F$ can be written as a graph
of $f$ for $t\in [0, \epsilon)$ and derive the equation satisfied
by $f$. It is given by a smooth family $F_t : \Sigma_1 \mapsto M$
which satisfies

\[(\frac{\partial F}{\partial t})^\perp
=H\] where $(\cdot)^\perp$ denotes the projection onto $N\Sigma$
and $H$ is the mean curvature vector of $F_t(\Sigma_1)=\Sigma_t$.
By the definition of the mean curvature vector, this equation is
equivalent to
\[(\frac{\partial F}{\partial t})^\perp
=(\Lambda^{ij} \nabla_{\frac{\partial F}
{\partial x^i}}^M\frac{\partial F}{\partial x^j})^
\perp\]
where $\Lambda^{ij}$ is the inverse to
 the induced metric $\Lambda_{ij}$
on $\Sigma$.
 \[\Lambda_{ij}=<\frac{\partial F}{\partial x^i},
 \frac{\partial F}{\partial x^j}>\]

When $\Sigma_1$ and $\Sigma_2$ are
both Euclidean space,  $\frac{\partial F}{\partial t}$ and
$\nabla_{\frac{\partial F}
{\partial x^i}}^M\frac{\partial F}{\partial x^j}$
are both in $T\Sigma_2$. Since the projection
to the normal part is an isomorphism when
restricted to $T\Sigma_2$,
\[\frac{\partial F}{\partial t}
=\Lambda^{ij} \nabla_{\frac{\partial F}
{\partial x^i}}^M\frac{\partial F}{\partial x^j}\]

If we write $F_t=id \times f_t$, then $\Lambda_{ij}=g_{ij}
+h_{\alpha\beta}\frac{\partial f^\alpha}{\partial x^i}
\frac{\partial f^\beta}{\partial x^j} =g_{ij}+(f^* h)_{ij}$. $f_t$
satisfies the following nonlinear partial differential equations.

\[\frac{\partial f^\alpha}{\partial t}
=( g_{ij}
+h_{\gamma\beta}\frac{\partial f^\gamma}{\partial x^i}
\frac{\partial f^\beta}{\partial x^j})^{-1}
\frac{\partial^2 f^\alpha}{\partial x^i \partial x^j}
\]

\section{Evolution equation for parallel form}
In this section, we calculate the evolution equation of the
restriction of a parallel $n$-form to an $n$-dimensional
submanifold moving by the mean curvature flow.

We assume $M $ is a Riemannian manifold with a parallel $n$ form
$\Omega$. Let $F:\Sigma\mapsto M$ be an isometric immersion of an
$n$-dimensional submanifold. We choose orthonormal frames
$\{e_i\}_{i=1\cdots n}$ for $T\Sigma$ and
$\{e_\alpha\}_{\alpha=n+1, \cdots, n+m}$ for $N\Sigma$. The
convention that $i, j, k, \cdots$ denote tangent indexes and
$\alpha, \beta, \gamma \cdots $ denote normal indexes is followed.

We first calculate the covariant derivative of
the restriction of $\Omega$ on $\Sigma$.

\[
\begin{split}
&(\nabla_{e_k}^\Sigma \Omega)(e_{i_1},\cdots,
 e_{i_n})\\
&=e_k(\Omega(e_{i_1}, \cdots, e_{i_n}))-
\Omega(\nabla_{e_k}^\Sigma e_{i_1}, \cdots, e_{i_n})-
\cdots-\Omega(e_{i_1}, \cdots, \nabla_{e_k}^\Sigma e_{i_n})\\
&=\Omega(\nabla_{e_k}^M e_{i_1}-
\nabla_{e_k}^\Sigma e_{i_1}, \cdots, e_{i_n})+\cdots
+\Omega(e_{i_1}, \cdots, \nabla_{e_k}^M e_{i_n}-
\nabla_{e_k}^\Sigma e_{i_n})
\end{split}
\]
where we have used $\nabla^M_{e_k}\Omega=0$ because $\Omega $ is
parallel. This equation can be abbreviated using the second
fundamental form of $F$, $h_{\alpha ij} =<\nabla_{e_i}^M e_j,
e_\alpha>$.

\begin{equation}\label{gradient}
\Omega_{i_1 \cdots i_n, k}
=\Omega_{\alpha i_2 \cdots i_n}h_{\alpha i_1 k}
+\cdots +\Omega_{i_1 \cdots i_{n-1} \alpha}
h_{\alpha i_n k}
\end{equation}

Likewise,
\begin{equation}\label{alphak}
\Omega_{\alpha i_2 \cdots i_n,k}
=-\Omega_{j i_2\cdots i_n}h_{\alpha jk}
+\Omega_{\alpha \beta i_3 \cdots i_n}
h_{\beta i_2 k}
+\cdots
+\Omega_{\alpha i_2 \cdots i_{n-1} \beta}
h_{\beta i_n k}
\end{equation}

Take the covariant derivative of
equation (\ref{gradient}) with respect
to $e_k$ again,

\begin{equation} \label{kk}
\begin{split}
\Omega_{i_1 \cdots i_n, kk}
&= \Omega_{\alpha i_2 \cdots i_n,k}h_{\alpha i_1 k}
+\cdots +\Omega_{i_1 \cdots i_{n-1} \alpha,k}
h_{\alpha i_n k}\\
&+\Omega_{\alpha i_2 \cdots i_n}
h_{\alpha i_1 k,k}
+\cdots +\Omega_{i_1 \cdots i_{n-1} \alpha}
h_{\alpha i_n k,k}
\end{split}
\end{equation}

Plug equation (\ref{alphak}) into (\ref{kk}) and apply the
Codazzi equation $h_{\alpha ki,k}=h_{\alpha ,i}+R_{\alpha kki}$
where $R$ is the curvature operator of $M$. Now we specialize to
$i_1=1, \cdots, i_n=n$.

\[
\begin{split}
(\Delta^\Sigma\Omega)_{1\cdots n} &=(-\Omega_{j 2\cdots
n}h_{\alpha jk} +\Omega_{\alpha \beta 3\cdots n}h_{\beta 2k}
+\cdots+ \Omega_{\alpha 2\cdots (n-1)\beta}
h_{\beta nk})h_{\alpha 1k}\\
&+\cdots
+(\Omega_{\beta 2\cdots (n-1)\alpha}h_{\beta 1k}
+\cdots +\Omega_{1\cdots (n-2)\beta \alpha}
h_{\beta (n-1) k}
-\Omega_{1\cdots (n-1) j}h_{\alpha jk})
h_{\alpha nk}\\
&+\Omega_{\alpha 2\cdots n}h_{\alpha,1}
+\cdots+\Omega_{1\cdots (n-1)\alpha}
h_{\alpha,n}\\
&+\Omega_{\alpha 2 \cdots n}R_{\alpha kk1}
+\cdots +\Omega_{1\cdots (n-1)\alpha}R_{\alpha
kkn}
\end{split}
\]
where $\Delta^\Sigma\Omega$ is the rough Laplacian, i.e.

\[(\Delta^\Sigma\Omega)_{1\cdots
n}=(\nabla^\Sigma_{e_k}\nabla^\Sigma_{e_k}\Omega)(e_1,\cdots,e_n)\]

Since $\Sigma$ is of dimension $n$, after
grouping terms we have
\begin{equation}\label{lap}
\begin{split}
(\Delta^\Sigma \Omega)_{1\cdots n} &=-\Omega_{1 2\cdots n}
\sum_{\alpha, k}(h_{\alpha 1k}^2 +\cdots+h_{\alpha nk}^2)\\
&+2\sum_{\alpha, \beta, k}[\Omega_{\alpha \beta 3\cdots n}
h_{\alpha 1k}h_{\beta 2k}
+ \Omega_{\alpha 2 \beta\cdots n}
h_{\alpha 1k}h_{\beta 3k}
+\cdots
+\Omega_{1\cdots (n-2) \alpha \beta}
h_{\alpha (n-1)k}
h_{\beta nk}]\\
&+\sum_{\alpha, k}\Omega_{\alpha 2\cdots n}h_{\alpha,1}
+\cdots+\Omega_{1\cdots (n-1)\alpha}
h_{\alpha,n}\\
&+\sum_{\alpha,k} \Omega_{\alpha 2 \cdots n}R_{\alpha kk1}
+\cdots +\Omega_{1\cdots (n-1)\alpha}R_{\alpha
kkn}
\end{split}
\end{equation}
We notice that $(\Delta^\Sigma \Omega)_{1\cdots n}=\Delta
(\Omega(e_1,\cdots, e_n))$, where the $\Delta$ on the right hand
side is the Laplacian of functions on $\Sigma$.

The terms in the bracket are formed in the
following way. Choose two different indexes
from $1$ to $n$, replace the smaller one by
$\alpha$ and the larger one by $\beta$.
There are a total of $\frac{n(n-1)}{2}$ such terms.

Now we consider the mean curvature flow of $\Sigma$ in $M$ by
$\frac{d}{dt} F_t= H_t$. In the following, we shall denote the
image of $F_t$ by $\Sigma_t$.  Notice that here we require the
velocity vector is in the normal direction. The evolution equation
of $\Omega_{1\cdots n}$ can be
 calculated as the following. We work in a
 local
coordinate $\{\partial_i
=\frac{\partial}{\partial x^i}\}$ on $\Sigma$, then

\[
\begin{split}
&\frac{d}{dt}\Omega(\partial_1,\cdots, \partial_n)\\
&=\Omega((\nabla_{\partial_1}H)^N,\partial_2,
\cdots, \partial_n)+\cdots+
\Omega(\partial_1,\partial_2,
\cdots, (\nabla_{\partial_n} H)^N)\\
&+\Omega((\nabla_{\partial_1}H)^T,\partial_2,
\cdots, \partial_n)+\cdots+
\Omega(\partial_1,\partial_2,
\cdots, (\nabla_{\partial_n} H)^T)
\end{split}
\]
Since $\frac{d}{dt} g_{ij}=
<(\nabla_{\partial_i}H)^T,\partial_j>$, if
we choose a orthonormal frame and evolve
the frame with respect to time so that it
remains orthonormal, the terms in the last
line vanish.

\[\frac{d}{dt}\Omega_{1 \cdots n}
=\Omega_{\alpha 2\cdots n} h_{\alpha,1}
+\cdots +\Omega_{1\cdots (n-1)\alpha}
h_{\alpha,n}
\]

 Combine this with equation (\ref{lap}) we
 get the parabolic equation
satisfied by $\Omega_{1\cdots n}$.
\begin{pro}
If $\Sigma_t$ is an $n$-dimensional mean curvature flow in $M$ and
$\Omega$ is a parallel $n$-form on $M$. Then $\Omega_{1,\cdots, n}
=\Omega(e_1, \cdots, e_n)$ satisfies

\begin{equation}\label{main}
\begin{split}
\frac{d}{dt}\Omega_{1 \cdots n}
&=\Delta\Omega_{1\cdots n}
+\Omega_{1 \cdots n}
(\sum_{\alpha ,i,k}h_{\alpha i k}^2)\\
&-2\sum_{\alpha, \beta, k}[\Omega_{\alpha \beta 3\cdots n}
h_{\alpha 1k}h_{\beta 2k}
+ \Omega_{\alpha 2 \beta\cdots n}
h_{\alpha 1k}h_{\beta 3k}
+\cdots
+\Omega_{1\cdots (n-2) \alpha \beta}
h_{\alpha (n-1)k}
h_{\beta nk}]\\
&-\sum_{\alpha, k}[\Omega_{\alpha 2 \cdots n}R_{\alpha kk1}
+\cdots +\Omega_{1\cdots (n-1)\alpha}R_{\alpha
kkn}]
\end{split}
\end{equation}
where $\Delta$ denotes the time-dependent Laplacian on $\Sigma_t$.
\end{pro}

When $M=\Sigma_1\times\Sigma_2$ is product, the volume form
$\Omega_i$ of  each $\Sigma_i$ is a parallel form on $M$. In fact,
all the discussions in this paper apply to any locally Riemannian
product manifold. At any point $p$ on $\Sigma$, choose a oriented
orthonormal basis $e_1, \cdots e_n$ for $T_p \Sigma$. Then
$\Omega_1(T\Sigma)= \Omega_1(e_1,\cdots
e_n)=\Omega_1(\pi_1(e_1),\cdots ,\pi_1(e_n))$ is the Jacobian of
the projection from $T\Sigma$ to $T\Sigma_1$. We shall use
$*\Omega_1$ to denote this function as $p$ varies along $\Sigma$,
here $*$ is simply the Hodge operator with respect to the induced
metric. By the implicit function theorem, $*\Omega_1>0$ near $p$
if and only if $\Sigma$ is locally a graph over $\Sigma_1$ near
$p$.

In the following, we assume $\Sigma_t $ is locally a graph over
$\Sigma_1$ so that $*\Omega_1>0$ on $\Sigma_t $ and find
orthonormal bases for the tangent and normal bundle of
 $\Sigma_t$ so that we can represent the terms
 in equation (\ref{main}) in a better form.

 First, we need a simple linear algebra lemma. Let
 $V=V_1\times V_2$ be a product of inner
  product spaces $V_1$ and $V_2$ of dimension
  $n$ and $m$ respectively.
 Let $D:V_1 \mapsto V_2 $ be a linear
 transformation.

\begin{lem}
There exist orthonormal bases $\{a_i\}_{i=1
\cdots n}$ for $V_1$ and $\{a_\alpha\}_{
\alpha=n+1\cdots n+m}$ for
$V_2$ such that $\lambda_{i\alpha}
=<Da_i, a_\alpha>\geq 0 $ is diagonal.
\end{lem}

In fact, this is the Singular Value Decomposition and a proof is
available in e.g. \cite{leon}.

It is understood that if $n<m$, then $\lambda_{i\alpha}=0$ for
$\alpha>n$ and if $m<n$, then $\lambda_{i\alpha}=0$ for $i>m$.

Now let $T$ be the graph of $D$,
 i.e. $T=V_1+D(V_1)$. $N$ denotes the orthogonal
 complement of $T$. Let $\pi_i :V_1\times V_2
 \mapsto V_i$ be the projection map.
 We notice that there are isomorphism
 $\pi_1|_T:T\mapsto V_1$ and $\pi_2|_N: N
 \mapsto V_2$.

 In the later application, $V_i=T_p \Sigma_i$ and
 $D$ is given by the $df|_p$, the differential
 of $f$ at the point $p$, where $f$ is a locally defined map whose
 graph represents $\Sigma$ near $p$.

Now $\{e_i= \frac{1}{\sqrt{1+\sum_{\beta}\lambda_{i\beta}^2}}
(a_i+\sum_\beta \lambda_{i\beta}a_\beta)\}$ forms an orthonormal
basis for $T$ and $\{e_\alpha=
 \frac{1}{\sqrt{1+\sum_j \lambda_{j\alpha}^2}}
(a_\alpha-\sum_j \lambda_{j\alpha}a_j)\}$ an orthonormal basis for
$N$. It is not hard to check that $|\pi_1(e_i)|=\frac{1}
{\sqrt{1+\sum_{\beta} \lambda^2_{i\beta}}}$ and

\begin{equation}\label{pi1}
\begin{split}
\pi_1(e_\alpha)&=-\sum_j\lambda_{j\alpha}
\pi_1(e_j)\\
\pi_2(e_i)&=\sum_i\lambda_{i\beta}
\pi_2(e_\beta)\\
\end{split}
\end{equation}

With these bases, we can calculate the terms in equation
(\ref{main}) for $\Omega=\Omega_1$.

We first calculate the term
\[-2\sum_{\alpha, \beta, k}[\Omega_{\alpha \beta 3\cdots n}
h_{\alpha 1k}h_{\beta 2k}
+ \Omega_{\alpha 2 \beta\cdots n}
h_{\alpha 1k}h_{\beta 3k}
+\cdots
+\Omega_{1\cdots (n-2) \alpha \beta}
h_{\alpha (n-1)k}
h_{\beta nk}]
\]

By equation (\ref{pi1}),
\[\Omega_1(e_\alpha, e_\beta, e_3,\cdots, e_n)
=\Omega_1(\pi_1(e_\alpha), \pi_1(e_\beta),
\pi_1(e_3), \cdots, \pi_1(e_n))
=(\lambda_{1\alpha }\lambda_{2\beta }
-\lambda_{2\alpha }\lambda_{1\beta })
*\Omega_1
\]

Therefore, this term is

\begin{equation}\label{quad}
-2[\sum_{\alpha, \beta, k, i<j} (\lambda_{i\alpha
}\lambda_{j\beta } -\lambda_{j\alpha }\lambda_{i\beta }) h_{\alpha
ik} h_{\beta jk}] *\Omega_1
\end{equation}

As for the curvature term,

\[\Omega_{\alpha ,2, \cdots, n}
R_{\alpha kk1}
=\Omega_1(e_\alpha, e_2, \cdots, e_n)
R(e_\alpha, e_k, e_k ,e_1)
\]

We assume $\Sigma_i $ is of constant curvature
$k_i$, therefore

\[
\begin{split}
&R(e_\alpha, e_k, e_k ,e_1)\\
&=R_1(\pi_1(e_\alpha), \pi_1(e_k), \pi_1(e_k)
,\pi_1(e_1))
+R_2(\pi_2(e_\alpha), \pi_2(e_k), \pi_2(e_k)
,\pi_2(e_1))\\
&=k_1[<\pi_1(e_\alpha), \pi_1(e_k)>
<\pi_1(e_k), \pi_1(e_1)>
-<\pi_1(e_\alpha), \pi_1(e_1)>
<\pi_1(e_k), \pi_1(e_k)>]\\
&+k_2[<\pi_2(e_\alpha), \pi_2(e_k)>
<\pi_2(e_k), \pi_2(e_1)>
-<\pi_2(e_\alpha), \pi_2(e_1)>
<\pi_2(e_k), \pi_2(e_k)>]
\end{split}
\]

Notice that
\[<X,Y>=<\pi_1(X), \pi_1(Y)>
+<\pi_2(X), \pi_2(Y)>\]
since $T\Sigma_1 \perp T\Sigma_2$.

Therefore the second term is
\[
\begin{split}
&k_2\sum_k[(-<\pi_1(e_\alpha), \pi_1(e_k)>)
(\delta_{1k}-<\pi_1(e_k), \pi_1(e_1)>)\\
&-(-<\pi_1(e_\alpha), \pi_1(e_1)>)
(1-<\pi_1(e_k), \pi_1(e_k)>)]\\
&=k_2\sum_k[<\pi_1(e_\alpha), \pi_1(e_1)>
-\delta_{1k}<\pi_1(e_\alpha), \pi_1(e_k)>\\
&+<\pi_1(e_\alpha), \pi_1(e_k)>
<\pi_1(e_k), \pi_1(e_1)>
-<\pi_1(e_\alpha), \pi_1(e_1)>
<\pi_1(e_k), \pi_1(e_k)>]
\end{split}
\]
Plug in $\pi_1(e_\alpha)=-\lambda_{j\alpha }
\pi_1(e_j)$.
We get

\[
\begin{split}
&\sum_{\alpha, k}\Omega_{\alpha 2\cdots n}R_{\alpha kk1}\\
&=\sum_{\alpha, j, k}\lambda_{1\alpha }
\lambda_{j\alpha }
\{k_1[<\pi_1(e_j), \pi_1(e_k)><\pi_1(e_k),
\pi_1(e_1)>-<\pi_1(e_j), \pi_1(e_1)>
|\pi_1(e_k)|^2]\\
&+k_2[<\pi_1(e_j), \pi_1(e_1)>
-\delta_{1k}<\pi_1(e_j), \pi_1(e_k)>\\
&+<\pi_1(e_j), \pi_1(e_k)>
<\pi_1(e_k), \pi_1(e_1)>
-<\pi_1(e_j), \pi_1(e_1)>
|\pi_1(e_k)|^2]\} *\Omega_1
\end{split}
\]

Similarly we can write down other
terms and the curvature term becomes
\[
\begin{split}
-\sum_{ \alpha,i,j, k} \lambda_{i\alpha }
\lambda_{j\alpha }&\{k_1[<\pi_1(e_j), \pi_1(e_k)><\pi_1(e_k),
\pi_1(e_i)>-<\pi_1(e_j), \pi_1(e_i)>
|\pi_1(e_k)|^2]\\
&+k_2[<\pi_1(e_j), \pi_1(e_i)>
-\delta_{ik}<\pi_1(e_j), \pi_1(e_k)>\\
&+<\pi_1(e_j), \pi_1(e_k)>
<\pi_1(e_k), \pi_1(e_i)>
-<\pi_1(e_j), \pi_1(e_i)>
|\pi_1(e_k)|^2]\} *\Omega_1
\end{split}
\]

Notice that $\sum_\alpha \lambda_{i\alpha}
\lambda_{j\alpha}\not=0$ only if $i=j$.
Let $\sum_\alpha \lambda_{i\alpha}
\lambda_{j\alpha}=\delta_{ij}\lambda_i^2$
with $\lambda_i\geq 0$.

We can rearrange the basis $\{a_\alpha\}$ so that
\[\lambda_{i\alpha }=\delta_{\alpha,
n+i}\lambda_i\]

If we represent $\Sigma_t$ locally as the graph of a locally
defined map $f_t$, then $\lambda_i$'s are in fact the eigenvalues
of $\sqrt{(df_t)^T df_t}$.

\[
\begin{split}
&\sum_{\alpha, i,k}\lambda_{i\alpha }^2
\{k_1[|\pi_1(e_i)|^2|\pi_1(e_k)|^2
-<\pi_1(e_i), \pi_1(e_k)>^2]\\
&+k_2[|\pi_1(e_i)|^2 |\pi_1(e_k)|^2
-<\pi_1(e_i), \pi_1(e_k)>^2
+\delta_{ik}<\pi_1(e_i), \pi_1(e_k)>
-|\pi_1(e_i)|^2]\} *\Omega_1
\end{split}
\]

By the choice of  $\{e_i\}$, $<\pi_1(e_i), \pi_1(e_k)>=0$ unless
$i=k$, and if we write $\sum_{i}|\pi_1 (e_i)|^2=|\pi_1|^2$, then
the term is

\[
\begin{split}
&\sum_{\alpha, i}\lambda_{i\alpha }^2
|\pi_1(e_i)|^2\{k_1
(|\pi_1|^2-|\pi_1(e_i)|^2)
+k_2  (|\pi_1|^2-
|\pi_1(e_i)|^2+1-n)\} *\Omega_1
\end{split}
\]

Now we shall write the equation for
$*\Omega$ in terms of $\lambda_i$ and
the second fundamental form. Notice
that $|\pi_1(e_i)|^2=\frac{1}{1+\lambda_i
^2}$.

\begin{pro}\label{equation}
Suppose $M=\Sigma_1\times \Sigma_2$ and $\Sigma_i$ is a Riemannian
manifold of constant curvature $k_i$, $i=1, 2$. The volume form of
$\Sigma_i$ is denoted by $\Omega_i$.  Let $F_0:\Sigma
\hookrightarrow M$ be an embedding such that $\Sigma$ is locally a
graph over $\Sigma_1$. If each $\Sigma_t$ is  locally a graph over
$\Sigma_1$ along the mean curvature flow of $F_0$ for $t\in
[0,\epsilon)$, then $*\Omega=*\Omega_1$ satisfies the following
equation.
\begin{equation}\label{eq}
\begin{split}
&\frac{d}{dt}*\Omega=\Delta *\Omega
+*\Omega\{\sum_{\alpha, i, k}
h_{\alpha ik}^2-2\sum_{k, i<j}
\lambda_{ i}\lambda_{ j}h_{n+i,ik}
h_{n+j, jk}
+2\sum_{k, i<j}\lambda_{i}\lambda_{j}
h_{n+j, ik} h_{n+i, jk}\\
&+\sum_{ i}\frac{\lambda_{ i}^2}{1+
\lambda_i^2}[k_1
(\sum_{j\not= i}\frac{1}{1+\lambda_j^2})
+k_2
(1-n+\sum_{j\not= i}\frac{1}{1+\lambda_j^2})]
\}
\end{split}
\end{equation}
where $\lambda_i's$ are the eigenvalues of $\sqrt{(df_t)^T df_t}$
and $f_t$ is a locally defined map whose graph represents
$\Sigma_t$ locally.
\end{pro}

It is understood that in case
$n>m$, we pretend $h_{n+i, jk}=0$ for
$m<i\leq n$.

If $\Sigma_t$ is indeed the graph of a map $f_t$, then
$\lambda_i$'s are the eigenvalues of $\sqrt{(df_t)^T df_t}$. As a
comparison with the harmonic heat flow, $*\Omega=\frac{1}
{\sqrt{\prod_{i=1}^n(1+\lambda_{i}^2)}}$ and the energy density of
$f$ is  $|df|^2=\sum_{i=1}^n \lambda_i^2$ . A lower bound of
$*\Omega$ implies an upper bound for $|df|^2$.

When $k_1=k_2=c$, the equation becomes.

\begin{equation}\label{same}
\begin{split}
&\frac{d}{dt}*\Omega=\Delta *\Omega
+*\Omega\{\sum_{\alpha, i, k}
h_{\alpha ik}^2-2\sum_{k, i<j}
\lambda_{ i}\lambda_{ j}h_{n+i,ik}
h_{n+j, jk}
+2\sum_{k, i<j}\lambda_{i}\lambda_{j}
h_{n+j, ik} h_{n+i, jk}\\
&+c\sum_{ i}\frac{\lambda_{ i}^2}{1+
\lambda_i^2}[
(\sum_{j\not= i}\frac{2}{1+\lambda_j^2})
+1-n]
\}
\end{split}
\end{equation}

\section{Long time existence and Convergence}
In this section, we consider the mean curvature flow of $\Sigma$
in $M=\Sigma_1\times \Sigma_2$ in the case when $\Sigma$ is the
graph of $f:\Sigma_1\mapsto \Sigma_2$. In particular, we prove the
long time existence and convergence, assuming that
$\det(g_{ij}+(f^*h)_{ij})$ is less than $2$ initially. In our
notation, $\det(g_{ij}+(f^*h)_{ij})= \det(\delta_{ij}+<f_*(a_i),
f_*(a_j)>_h)$, where $\{a_i\}$ is any orthonormal basis for
$(\Sigma_1, g)$.

Let us explain the hypothesis $\det(g_{ij}+(f^*h)_{ij})<2$ in more
detail when $\Sigma_1$ and $\Sigma_2$ are both two-dimensional
surfaces. As remarked in \S 1, this condition is equivalent to the
Jacobian $J_1$ of the projection from $\Sigma$ onto $\Sigma_1$ is
greater than $\frac{1}{\sqrt{2}}$ and is slightly stronger than
$J_1>|J_2|$.

By the Singular Value Decomposition at any point $p\in \Sigma_1$,
we can choose orthonormal bases $\{a_1, a_2\}$ for $T_p\Sigma_1$
and $\{a_3, a_4\}$ for $T_{f(p)}\Sigma_2$ such that
$df|_p(a_1)=\lambda_1 a_3$ and $df|_p(a_2)=\lambda_2 a_4$. Then
$\det(g_{ij}+(f^*h)_{ij})=(1+\lambda_1^2)(1+\lambda_2^2)$. Let
$\Omega_1$ and $\Omega_2$ be the volume form of $\Sigma_1$ and
$\Sigma_2$ respectively. They can be extended as parallel forms on
$\Sigma_1\times\Sigma_2$. We also have the projections
$\pi_1:\Sigma_1\times\Sigma_2\mapsto \Sigma_1$ and
$\pi_2:\Sigma_1\times\Sigma_2\mapsto \Sigma_2$. At any point $(p,
f(p))\in \Sigma$ and any orthonormal basis $\{e_1, e_2\}$ for
$T_{(p,f(p))} \Sigma$, $\Omega_1(e_1, e_2)$ is the Jacobian of
$\pi_1|_\Sigma$, the restriction of $\pi_1$ to  $\Sigma$, and
$\Omega_2(e_1, e_2)$ is the Jacobian of $\pi_2|_\Sigma$. Now we
can take the orthonormal basis for $T_{(p, f(p))}\Sigma$ to
consist of
\[ e_1=\frac{1}{\sqrt{1+\lambda_1^2}}(a_1+\lambda_1 a_3),\,
e_2=\frac{1}{\sqrt{1+\lambda_2^2}}(a_2+\lambda_2 a_4)\]

Since $\Omega_1=a_1^*\wedge a_2^*$ and $\Omega_2=a_3^*\wedge
a_4^*$, we have $\Omega_1(e_1,
e_2)=\frac{1}{\sqrt{(1+\lambda_1^2)(1+\lambda_2^2)}}$ and
$\Omega_2(e_1, e_2)=\frac{\lambda_1 \lambda_2
}{\sqrt{(1+\lambda_1^2)(1+\lambda_2^2)}}$. The assumption
$(1+\lambda_1^2)(1+\lambda_2^2)<2$ is equivalent to

\begin{equation}\label{dom}
\Omega_1(e_1, e_2)>\frac{1}{\sqrt{2}}
\end{equation}

Taking into account of the fact that
$(\frac{1}{\sqrt{(1+\lambda_1^2)(1+\lambda_2^2)}})^2+(\frac{\lambda_1
\lambda_2 }{\sqrt{(1+\lambda_1^2)(1+\lambda_2^2)}})^2\leq 1$,
(\ref{dom}) implies the weaker condition
 \[\Omega_1(e_1,
e_2)>|\Omega_2(e_1, e_2)|\]

In fact the condition $ \Omega_1(e_1, e_2)>|\Omega_2(e_1, e_2)|$
is equivalent to $\Sigma$ being symplectic with respect to both
symplectic forms $\Omega_1+\Omega_2$ and $\Omega_1-\Omega_2$. This
is exactly the assumption in \cite{ke} where we prove the global
existence and convergence in the two-dimensional case. First we
prove Theorem A.

\vskip 10pt \noindent {\bf Theorem A.} {\it Let $(\Sigma_1, g)$
and $(\Sigma_2, h)$ be Riemannian manifolds of constant curvature
$k_1$ and $k_2$ respectively
 and $f$ be
a smooth map from $\Sigma_1$ to $\Sigma_2 $. Suppose $k_1\geq
|k_2|$. If $\,\det(g_{ij}+(f^*h)_{ij})<2$, the mean curvature flow
of the graph of $f$
 remains a graph and exists for all time.} \vskip 10pt

\begin{proof} Following the notation in the previous section
with $\Omega=\Omega_1$. It is not hard to see that
$*\Omega=\frac{1}{\sqrt{\det(g_{ij}+(f^*h)_{ij})}}=\frac{1}
{\sqrt{\prod_{i=1}^n(1+\lambda_{i}^2)}}$ and the assumption
implies $\prod_{i=1}^n(1+\lambda_{i}^2)\leq 2-\delta$ for some
$\delta>0$, and in particular $\sum_{i=1}^n \lambda_i^2 \leq
1-\delta$.

Now we take a look at the quadratic terms of the second
fundamental form in equation (\ref{eq}). First we assume $n\leq
m$, so $n+m\geq 2n$. The tangent indices $i, j ,k $ run from $1$
to $n$ and the normal index $\alpha $ runs from $n+1$ to $n+m$
unless they  are specified otherwise. We divide $\sum h^2_{\alpha
ik}$ into two parts:

\[\sum_{\alpha, i,k} h_{\alpha i k}^2=\sum_{n+1\leq \alpha \leq
2n, i,k} h_{\alpha ik}^2 +\sum_{2n<\alpha, i,k} h_{\alpha ik}^2\]

In the first summand, write $\alpha =n+j$ then $j$ runs from $1$
to $n$, therefore

\[\sum_{n+1\leq \alpha \leq 2n, i,k} h_{\alpha
ik}^2=\sum_{j,i,k} h^2_{n+j, i, k}=\sum_{i<j, k}
(h^2_{n+i,j,k}+h^2_{n+j,i,k})+\sum_{i,k} h^2_{n+i,i,k}\]

The quadratic terms of the second fundamental form in equation
(\ref{eq}) become.
\[
\begin{split}
&\sum_{\alpha, i, k}h_{\alpha ik}
^2-2\sum_{i<j, k}
\lambda_{ i}\lambda_{ j}h_{n+i,ik}
h_{n+j,jk}
+2\sum_{i<j, k}\lambda_{j}\lambda_{i}
h_{n+j, ik} h_{n+i, jk}\\
&=\delta |A|^2+(1-\delta)\sum_{\alpha >2n, i, k}h_{\alpha, ik}
^2+(1-\delta)\sum_{ i, k}h_{n+i, ik} ^2+(1-\delta)\sum_{i<j, k}
(h_{n+i, jk}^2
+ h_{n+j, ik}^2)\\
&-2\sum_{i<j, k}
\lambda_{ i}\lambda_{ j}h_{n+i,ik}
h_{n+j,jk}
+2\sum_{i<j, k}\lambda_{j}\lambda_{i}
h_{n+j, ik} h_{n+i, jk}\\
&\geq \delta |A|^2+(1-\delta)\sum_{\alpha >2n, i, k}h_{\alpha, ik}
^2+(\sum_{i,k} h_{n+i, ik}^2) (\sum_i
\lambda_i^2)+(1-\delta)\sum_{i<j, k} (h_{n+i, j,k}^2+h_{n+j, ik}^2)\\
&-2\sum_{i<j, k} \lambda_{ i}\lambda_j h_{n+i,ik} h_{n+j,jk}
-2(1-\delta)\sum_{i<j,k} |h_{n+j, ik} h_{n+i, jk}|
\end{split}
\]
where we have used $\sum_i \lambda_i^2 \leq 1-\delta $ and
$|\lambda_i \lambda_j|\leq 1-\delta$.

Drop the non-negative term $(1-\delta)\sum_{\alpha >2n, i,
k}h_{\alpha, ik} ^2$ and the last expression is no less than

\[
\begin{split}
&\delta|A|^2+(\sum_i \lambda_i h_{n+i, ik})^2 -2\sum_{i<j,
k}\lambda_i\lambda_j h_{n+i,ik}h_{n+j, jk}
+(1-\delta)\sum_{i<j,k}(|h_{n+i, jk}|-|h_{n+j, ik}|)^2\\
&\geq \delta |A|^2+\sum_i \lambda_i^2  h_{n+i, ik}^2
+(1-\delta)\sum_{i<j, k}(|h_{n+i, jk}|-|h_{n+j, ik}|)^2
\end{split}
\]
which is non-negative.

If $n>m$, since $h_{n+i,jk}=0$ for
$m<i\leq n$, the quadratic terms become
\[
\begin{split}
&\sum_{\alpha, i, k}h_{\alpha ik}
^2-2\sum_{1\leq i<j\leq m, k}
\lambda_{ i}\lambda_{ j}h_{n+i,ik}
h_{n+j,jk}
+2\sum_{1\leq i<j\leq m , k}\lambda_{j}\lambda_{i}
h_{n+j, ik} h_{n+i, jk}\\
&=\sum_{\alpha, m<i\leq n , k}h_{\alpha, ik}
^2+\sum_{ \alpha, 1\leq i\leq m, k }
h_{\alpha, ik}^2\\
&-2\sum_{1\leq i<j\leq m, k}
\lambda_{ i}\lambda_{ j}h_{n+i,ik}
h_{n+j,jk}
+2\sum_{1\leq i<j\leq m , k}\lambda_{j}\lambda_{i}
h_{n+j, ik} h_{n+i, jk}\\
&=\sum_{\alpha, m<i\leq n , k}
h^2_{\alpha, ik}+
\sum_{1\leq i\leq m, k} h_{n+i, i k}^2
+\sum_{1\leq i<j\leq m, k} h_{n+i, jk}^2
+\sum_{1\leq i<j\leq m, k} h_{n+j, ik}^2\\
&-2\sum_{1\leq i<j\leq m, k}
\lambda_{ i}\lambda_{ j}h_{n+i,ik}
h_{n+j,jk}
+2\sum_{1\leq i<j\leq m , k}\lambda_{j}\lambda_{i}
h_{n+j, ik} h_{n+i, jk}\\
\end{split}
\]
By a similar argument, this term is non-negative and bounded below
by $\delta |A|^2$ .

As for the curvature term, for each $i$
we have
\begin{equation}\label{positive}
\begin{split}
&k_1
(\sum_{j\not= i}\frac{1}{1+\lambda_j^2})
+k_2
(1-n+\sum_{j\not= i}\frac{1}{1+\lambda_j^2})\\
&=(k_1+k_2) (\sum_{j\not= i}\frac{1}{1+\lambda_j^2})
+k_2(1-n)\\
&\geq \frac{(k_1-k_2)}{2}(n-1)+(k_1+k_2)[
(\sum_{j\not= i}\frac{1}{1+\lambda_j^2})-\frac
{n-1}{2}]
\end{split}
\end{equation}
Because each $\lambda_j^2$ is less than $1$ and $k_1\geq |k_2|$,
this term is nonnegative. When $k_1+k_2>0$, this term is indeed
strictly positive.

By Proposition \ref{equation} and the previous
paragraph,
$*\Omega$ satisfies the differential inequality.

\begin{equation}\label{omega}
\frac{d}{dt}*\Omega\geq \Delta *\Omega +\delta |A|^2
\end{equation}
According to the maximum principle for parabolic equations,
$\min_{\Sigma_t} *\Omega$ is nondecreasing in time. In particular,
$*\Omega$ has a positive lower bound. Since $*\Omega$ is the
Jacobian of the projection map from $\Sigma_t$ to $\Sigma_1$, by
the implicit function theorem, this implies $\Sigma_t$ remains the
graph of a map $f_t:\Sigma_1\mapsto \Sigma_2$ whenever the flow
exists.

Now we isometrically embed $M=\Sigma_1\times\Sigma_2$ into $\R^N$.
The mean curvature flow equation in terms of the coordinate
function $F(x,t)$ in $\R^N$ becomes
\[ \frac{d}{dt} F(x,t)=H=\overline{H}+E\]
where $H\in TM\slash T\Sigma$ is the mean curvature vector of
$\Sigma_t$ in $M$ and $\overline{H}\in T\R^N\slash T\Sigma$ is the
mean curvature vector of $\Sigma_t$ in $\R^N$.

To detect a possible singularity at $(y_0, t_0)$, recall the
(n-dimensional) backward heat kernel $\rho_{y_0, t_0}$ at $(y_0,
t_0)$ introduced by Huisken \cite{hu2}.

\[
\rho_{y_0, t_0}(y,t)=\frac{1}{(4\pi(t_0-t))^{n\over 2}} \exp
(\frac{-|y-y_0|^2}{4(t_0-t)})
\]
The monotonicity formula of Huisken asserts $\lim_{t\rightarrow
t_0}\int\rho_{y_0, t_0} d\mu_t$ exists. $\rho_{y_0,t_0}$ satisfies
the following backward heat equation derived in \cite{ke} along
the mean curvature flow. Here $\nabla$ and $\Delta$ represent the
covariant derivative and the Laplace operator of the induced
metric on $\Sigma_t$ respectively.

\begin{equation}\label{bh}
\begin{split}
\frac{d}{dt}\rho_{y_0,t_0}&=- \Delta\rho_{y_0,t_0}
-\rho_{y_0,t_0}( \frac{|F^\perp|^2}{4(t_0-t)^2} +\frac{F^\perp
\cdot \overline{H}}{t_0-t} +\frac{F^\perp \cdot E}{2(t_0-t)})
\end{split}
\end{equation}
where $F^\perp$ is the component of $F\in T\R^N$ in $T\R^N \slash
T\Sigma_t$.

Recall that

\[
\frac{d}{dt}d\mu_t =-|H|^2 d\mu_t =-\overline{H}\cdot
(\overline{H}+E) d\mu_t
\]

Combine this equation with equations (\ref{omega}) and (\ref{bh}),
we get

\begin{equation}
\begin{split}
&\frac{d}{dt}\int (1-*\Omega)
\rho_{y_0,t_0}\, d\mu_t\\
\leq &\int[\Delta
(1-*\Omega)-\delta |A|^2]\rho_{y_0,t_0} \,d\mu_t \\
&-\int (1-*\Omega)[\Delta\rho_{y_0,t_0} +\rho_{y_0,t_0}(
\frac{|F^\perp|^2}{4(t_0-t)^2} +\frac{F^\perp \cdot
\overline{H}}{t_0-t} +\frac{F^\perp \cdot E}{2(t_0-t)})]d\mu_t\\
&-\int(1-*\Omega)[\overline{H}(\overline{H}+E)]\rho_{y_0,
t_0}d\mu_t
\end{split}
\end{equation}

By rearranging terms, the right hand side can be written as
\[
\begin{split}
&\int[\Delta(1-*\Omega)\rho_{y_0,t_0}-(1-*\Omega) \Delta
\rho_{y_0, t_0}] \,d\mu_t
-\delta\int |A|^2\rho_{y_0, t_0} d\mu_t \\
&-\int (1-*\Omega)\rho_{y_0,t_0}[ \frac{|F^\perp|^2}{4(t_0-t)^2}
+\frac{F^\perp \cdot \overline{H}}{t_0-t} +\frac{F^\perp \cdot
E}{2(t_0-t)}+|\overline{H}|^2+\overline{H}\cdot E]d\mu_t
\end{split}
\]

The first term vanishes by integration by parts and the third term
can be completed square. Therefore
\[
\begin{split}
&\frac{d}{dt}\int (1-*\Omega)
\rho_{y_0,t_0}\, d\mu_t\\
\leq &-\delta\int |A|^2\rho_{y_0, t_0} d\mu_t-\int
(1-*\Omega)\rho_{y_0,t_0}\left|\frac{|F^\perp|}{2(t_0-t)}
+\overline{H} +\frac{ E}{2}\right|^2 d\mu_t\\
&+ \int(1-*\Omega)\rho_{y_0, t_0} \left|\frac{E}{2}\right|^2
d\mu_t
\end{split}
\]

Since $E$ is bounded and $\int(1-*\Omega)\rho_{y_0, t_0}
d\mu_t\leq \int\rho_{y_0, t_0} d\mu_t$ is finite, this implies
\[
\begin{split}
&\frac{d}{dt}\int (1-*\Omega)
\rho_{y_0,t_0}\, d\mu_t\\
\leq &C-\delta\int |A|^2\rho_{y_0, t_0} d\mu_t
\end{split}
\]
for some constant $C$. From this we see that $\lim_{t\rightarrow
t_0} \int (1-*\Omega)\rho_{y_0, t_0} d\mu_t$ exists.

For $\lambda>1$, the parabolic dilation $D_\lambda$ at $(y_0,
t_0)$ is defined by

\begin{equation}
\begin{split}
D_\lambda :\,\R^N \times[0, t_0) &\rightarrow
\R^N \times [-\lambda^2 t_0, 0)\\
(y,t)&\rightarrow (\lambda(y-y_0), \lambda^2 (t-t_0))
\end{split}
\end{equation}

Let $\mathcal{S}\subset \R^N\times[0, t_0)$ be the total space of
the mean curvature flow , we shall study  the flow
$\mathcal{S}^\lambda=D^\lambda (\mathcal{S})\subset \R^N\times
[-\lambda^2 t_0, 0)$ . Denote the new time parameter by $s$, then
$t=t_0+\frac{s}{\lambda^2}$. Let $d\mu^\lambda_s$ denote the
induced volume form on $\Sigma$ by $F^\lambda_s=\lambda
F_{t_0+\frac{s}{\lambda^2}}$. The image of $F^\lambda_s$ is the
$s$-slice of $\mathcal{S}^\lambda$ and is denoted by
$\Sigma^\lambda_s$. Therefore,
\[
\begin{split}
&\frac{d}{ds}\int
 (1-*\Omega) \rho_{0,0}\,
d\mu_{s}^{\lambda}\\
=&\frac{1}{\lambda^2}\frac{d}{dt}\int  (1-*\Omega)
\rho_{y_0,t_0}\,
d\mu_{t}\\
\leq &\, \frac{C}{\lambda^2}-\frac{\delta}{\lambda^2} \int
\rho_{y_0, t_0} |A|^2
\,d\mu_t\\
\end{split}
\]
Notice that $*\Omega $ is a invariant under the parabolic
dilation. It is not hard to check that
\[
\frac{1} {\lambda^2} \int \rho_{y_0, t_0} |A|^2\,d\mu_t =\int
\rho_{0,0}|A|^2 \,d\mu_{s}^{\lambda}
\]
This is because $\rho_{y_0, t_0} d\mu_t$ is invariant under the
parabolic scaling and the norm of second fundamental form scales
like the inverse of the distance.

Therefore

\[
\begin{split}
&\frac{d}{ds}\int
 (1-*\Omega) \rho_{0,0}\,
d\mu_{s}^{\lambda}\\
\leq &\frac{C}{\lambda^2}-\delta \int \rho_{0,0} |A|^2\,
d\mu_s^\lambda
\end{split}
\]

This reflects the correct scaling for the parabolic blow-up.

Take any $\tau>0$ and integrate from $-1-\tau$ to $-1$.

\begin{equation}\label{delta}
\begin{split}
&\delta\int_{-1-\tau}^{-1} \int \rho_{0,0}|A|^2
\,d\mu_{s}^{\lambda}ds\\
&\leq \int
 (1-*\Omega) \rho_{0,0}\,
d\mu_{-1}^{\lambda}-\int
 (1-*\Omega) \rho_{0,0}\,
d\mu_{-1-\tau}^{\lambda} +\frac{C}{\lambda^2}
\end{split}
\end{equation}
Notice that

\[\int  (1-*\Omega)
\rho_{0, 0} d\mu_{s} ^{\lambda} =\int (1-*\Omega)\rho_{y_0, t_0}
d\mu_{t_0+\frac{s}{\lambda^2}}
\]

By the fact that $\lim_{t\rightarrow t_0} \int
(1-*\Omega)\rho_{y_0, t_0} d\mu_t$ exists, the right hand side of
equation (\ref{delta}) approaches zero as $\lambda \rightarrow
\infty$. Take a sequence $\lambda_j \rightarrow \infty$, for a
fixed $\tau
>0$,

\[\int_{-1-\tau}^{-1} \int
\rho_{0, 0} |A|^2 d\mu_s^{\lambda_j}ds\leq C(j)
\]
where $C(j)\rightarrow 0$ as $\lambda_j \rightarrow \infty$.

Choose $\tau_j \rightarrow 0$ such that $\frac{C(j)}{\tau_j}
\rightarrow 0$ and $s_j \in [-1-\tau_j, -1]$ so that

\[\int \rho_{0, 0} |A|^2 d\mu_{s_j}^{\lambda_j}
\leq \frac{C(j)}{\tau_j}
\]

We investigate this inequality more carefully. Notice that

\[\rho_{0,0}(F^{\lambda_j}_{s_j})=\frac{1}{4\pi
(-s_j)}\exp(\frac{-|F^{\lambda_j}_{s_j}|^2}{4(-s_j)})\] where
$F^{\lambda_j}_{s_j}=\lambda_j F_{t_0+\frac{s_j}{\lambda_j^2}}$.

If we consider for any $R>0$, the ball of radius $R$,
$B_R(0)\subset \R^N$, when $j$ is large enough, we may assume
$-1<s_j<-\frac{1}{2}$, then

\[\int \rho_{0,0} |A|^2 d\mu_{s_j}^{\lambda_j}
\geq
\frac{1}{2\pi}\exp(\frac{-R^2}{2})\int_{\Sigma_{s_j}^{\lambda_j}\cap
B_R(0)}|A|^2 d\mu^{\lambda_j}_{s_j}\]

This implies for any compact set $K\subset \R^N$,

\begin{equation}\label{a2}
\int_{\Sigma_{s_j}^{\lambda_j}\cap K}|A|^2
d\mu_{s_j}^{\lambda_j}\rightarrow 0 \,\text{ as}\, j\rightarrow
\infty
\end{equation}

Now we claim in the rest of the proof this together with the fact
that $*\Omega $ has a positive lower bound imply  $
\lim_{j\rightarrow \infty}\int \rho_{y_0,
t_0}d\mu_{t_0+\frac{s_j}{\lambda_j^2}} =\lim_{j\rightarrow
\infty}\int \rho_{0, 0}d\mu^{\lambda_j}_{s_j}\leq 1$. We may
assume the origin is a limit point of $\Sigma_{s_j}^{\lambda_j}$,
otherwise the limit is zero and there is nothing to be proved.

$*\Omega$ is in fact the Jacobian of the projection
$\pi_1:\Sigma_t \mapsto \Sigma_1$. Each $\Sigma_t$ can be written
as the graph of a map $f_t :\Sigma_1 \mapsto \Sigma_2$ with
uniformly bounded $|d f_t|$.  This is because
$\det(g_{ij}+(f_t^*h)_{ij})=\prod_{i=1}^n (1+\lambda_i^2) $ is
bounded and  $\prod_{i=1}^n (1+\lambda_i^2) \geq 1+\sum_{i=1}^n
\lambda_i^2=1+|d f_t|^2$. Denote $f_{t_0+\frac{s_j}{\lambda_j^2}}$
by $f_j$. Now we consider the blow up of the graph of $f_j$ in
$\R^N$ by $\lambda_j$. This is the graph of the function
$\widetilde{f}_j$ defined on $\lambda_j \Sigma_1\subset \R^N$
which corresponds to a part of $\Sigma_{s_j}^{\lambda_j}$. Now
$|d\widetilde{f}_j|$ is also uniformly bounded and our assumption
on $\Sigma_{s_j}^{\lambda_j}$ implies $\lim_{j\rightarrow \infty}
\widetilde{f}_j(0)=0$. Therefore we may assume $\widetilde{f}_j
\rightarrow \widetilde{f}_\infty$ in $C^\alpha$ on compact sets.
$\widetilde{f}_\infty$ is an entire graph defined on $\R^n$.

Other the other hand,

\begin{equation}\label{A}
|A|_j\leq |\nabla d \widetilde{f}_j|\leq (1+|d
\widetilde{f}_j|^2)^{3\over 2}|A|_j
\end{equation}
where $|A|_j$ is the norm of the second fundamental form of
$\Sigma_{s_j}^{\lambda_j}$ and $|\nabla d\widetilde{f}_j|$ is the
norm of the covariant derivatives of $d\widetilde{f}_j$.
Inequalities (\ref{A}) can be derived as equations (29) on page 31
of \cite{il}.

Use the equation (\ref{a2}), we can show $ \widetilde{f}_j
\rightarrow
 \widetilde{f}_\infty$ in
$C^\alpha\cap W^{1,2}_{loc}$ and $\widetilde{f}_\infty$ has
vanishing second derivatives. This implies
$\Sigma^{\lambda_j}_{s_j} \rightarrow \Sigma^\infty_{-1}$ as Radon
measures and $\Sigma^\infty_{-1}$ is the graph of a linear
function. Therefore

\[\lim_{j\rightarrow \infty}\int \rho_{0,0}
d\mu_{s_j}^{\lambda_j} =\int \rho_{0,0} d\mu^\infty_{-1}=1\] This
implies \[\lim_{j\rightarrow \infty}\int \rho_{y_0, t_0}
d\mu_{t_0+\frac{s_j}{\lambda_j^2}} =\lim_{t\rightarrow t_0}\int
\rho_{y_0, t_0} d\mu_t=1\]

The regularity now follows from White's theorem \cite{w} which
asserts $(y_0, t_0)$ is a regular point whenever
$\lim_{t\rightarrow t_0}\int \rho_{y_0, t_0} d\mu_t\leq 1+\epsilon
$.

\end{proof}

\vskip 10pt \noindent {\bf Theorem B.} {\it  Suppose $k_1\geq
|k_2|$ and $k_1+k_2
>0$. If $\det(g_{ij}+(f^*h)_{ij})<2$,
then the flow exists for all time and the corresponding map
converges to a constant map at infinity.} \vskip 10pt

\begin{proof} Long time existence is already proved in Theorem A. Since
$*\Omega=\frac{1}{\sqrt{\det(g_{ij}+(f^*h)_{ij})}}=\frac{1}
{\sqrt{\prod_{i=1}^n(1+\lambda_{i}^2)}}$, the assumption is
equivalent to $*\Omega>\frac{1}{\sqrt{2}}$.

By equation (\ref{positive}), we have

\[k_1 (\sum_{j\not=
i}\frac{1}{1+\lambda_j^2}) +k_2 (1-n+\sum_{j\not=
i}\frac{1}{1+\lambda_j^2}) \geq 2c_1\]
for any $i$, where $c_1$ is
a constant that depends on the initial condition. By Proposition
\ref{equation}, $*\Omega$ satisfies

\[\frac{d}{dt}*\Omega\geq \Delta
*\Omega+2c_1\sum_{i=1}^n \frac{\lambda_i^2} {1+\lambda_i^2}\]

That $\lambda_i^2<1$ for each $i$ implies

\[\frac{d}{dt}*\Omega\geq \Delta
*\Omega+c_1\sum_{i=1}^n \lambda_i^2\]

Since $\lambda_i^2<1$
for each $i$, we have
\[1+c_2\sum_{i=1}^n \lambda_i^2\geq
\prod_{i=1}^n (1+\lambda_i^2)\geq 1+\sum_{i=1}^n \lambda_i^2\]
where $c_2$ is a constant that depends on $n$. Therefore
$\sum_{i=1}^n \lambda_i^2 \geq
\frac{1}{c_2}(\frac{1}{(*\Omega)^2}-1)$. and $*\Omega$ satisfies
\[\frac{d}{dt}*\Omega\geq \Delta
*\Omega+c_3(\frac{1}{(*\Omega)^2}-1)\]

By the comparison theorem for parabolic equations,
$\min_{\Sigma_t} *\Omega$ is non-decreasing in $t$ and
 $\min_{\Sigma_t} *\Omega \rightarrow 1$ as $t\rightarrow \infty$.

  To prove convergence at infinity, we first show that $\max_{x\in
\Sigma_t}|A|^2(x) \rightarrow 0$ as $t\rightarrow \infty$. We need
to take a look at the quadratic term of the second fundamental
form in equation (\ref{eq}).

 Let $\Lambda=(\lambda_{i \alpha})$ be a matrix and

\[Q(x)=\sum_{i, \alpha } x_{i \alpha }^2
-2\sum_{\alpha, \beta, i<j}
(\lambda_{i\alpha}\lambda_{j\beta }
-\lambda_{j\alpha }\lambda_{i\beta })
x_{i\alpha }x_{j\beta }
\]
be the quadratic form defined for
$x=(x_{i\alpha })\in \R^n \times \R^m$.
Choose $\epsilon$ small enough such that
$Q(x)>\frac{1}{2} |x|^2$
when $|\Lambda|^2\leq \epsilon$.

The quadratic term of the second fundamental form in equation
(\ref{eq}) in the original index (see also equation (\ref{quad}))
is

\begin{equation}
[\sum_{\alpha, i, k} h_{\alpha ik}^2-2\sum_{\alpha, \beta, i<j, k}
(\lambda_{i\alpha }\lambda_{j\beta } -\lambda_{j\alpha
}\lambda_{i\beta }) h_{\alpha ik} h_{\beta jk}] *\Omega_1
\end{equation}

Now take $\epsilon_2<\epsilon$. There exists
a time $T$ such that $*\Omega>
\frac{1}{\sqrt{1+\epsilon_2}}$
and $\sum \lambda_i^2 <\epsilon_2$ for $t>T$.
Therefore we have

\[\frac{d}{dt}*\Omega\geq \Delta *\Omega
+\frac{1}{2} *\Omega |A|^2\]

Let $\eta=*\Omega$, then by equation (\ref{gradient})
\[
\begin{split}
|\nabla \eta|^2&=\sum_k(\sum_\alpha
\Omega_{\alpha 2\cdots n}h_{\alpha 1k}
+\cdots+\Omega_{1\cdots, n-1\alpha}
h_{\alpha nk})^2\\
&\leq n \sum_k[\sum_\alpha
(\Omega_{\alpha 2\cdots n}h_{\alpha 1k})^2
+\cdots+(\Omega_{1\cdots, n-1\alpha}
h_{\alpha nk})^2]\\
&\leq n\sum_k[\sum_\alpha
(\lambda_{\alpha 1} h_{\alpha 1k})^2
+\cdots+(\lambda_{\alpha n}
h_{\alpha nk})^2](*\Omega)^2\\
\end{split}
\]
Therefore
\begin{equation}\label{geta}
|\nabla \eta|^2 \leq n\,\epsilon_2\,\eta^2 |A|^2
\end{equation}

 Let $p$ be a positive number to be determined, $\eta^p$
satisfies

\[
\begin{split}
\frac{d}{dt}\eta^p&=p\eta^{p-1}\frac{d}{dt}\eta\\
&\geq p\eta^{p-1}(\Delta \eta+\frac{1}{2}\eta
|A|^2)\\
&=\Delta \eta^p-p(p-1)\eta^{p-2}|\nabla \eta|^2
+\frac{p}{2}\eta^p|A|^2\\
\end{split}
 \]
Use the inequality (\ref{geta}), we get
\[
\begin{split}
\frac{d}{dt}\eta^p\geq \Delta \eta^p +[\frac{p}{2}-p(p-1)
n\epsilon_2 ]\eta^p|A|^2
\end{split}
 \]

Recall from \cite{ke} that $|A|^2$ satisfies the following
equation along the mean curvature flow.

\begin{equation}\label{|A|^2}
\begin{split}
\frac{d}{dt}|A|^2
&=\Delta |A|^2 -2|\nabla A|^2
+[(\nabla_{\pk}^M R)_{\alpha ijk}
+(\nabla_{\pj}^M R)_{\alpha kik}]h_{\alpha ij}\\
&-2R_{lijk}h_{\alpha lk}h_{\alpha ij}
+4R_{\alpha \beta jk}h_{\beta ik}h_{\alpha ij}
-2R_{lkik}h_{\alpha lj}h_{\alpha ij}
+R_{\alpha k\beta k}h_{\beta ij}h_{\alpha ij}\\
&+\sum_{\alpha,\gamma, i,m}
(\sum_k h_{\alpha ik}h_{\gamma mk}
-h_{\alpha mk}h_{\gamma ik})^2
+\sum_{i,j,m,k}(\sum_{\alpha} h_{\alpha ij}
h_{\alpha mk})^2
\end{split}
\end{equation}
where $R_{ABCD}$ is the curvature tensor
and $\nabla^M$ is the covariant derivative of
$M$.

In our case, the curvature operator is
parallel and it follows that $|A|^2$ satisfies

\[\frac{d}{dt}|A|^2\leq \Delta |A|^2
 -2|\nabla A|^2+K_1|A|^4+K_2 |A|^2\]
for $K_1$ , $K_2$ constants that depend on
the dimensions of $\Sigma_1$ and $\Sigma_2$.
Applying the same technique in \cite{lg}
to calculate $\eta^{-2p}|A|^2$ we get

\[
\begin{split}
\frac{d}{dt}(\eta^{-2p}|A|^2)
&\leq \Delta (\eta^{-2p}|A|^2)-\eta^{2p}
 \nabla(\eta^{-2p})\cdot
\nabla(\eta^{-2p}|A|^2) \\
&+\eta^{-2p}
[K_1|A|^4+K_2|A|^2-2|A|^4
(\frac{p}{2}-p(p-1)n \epsilon_2) ]
\end{split}
\]
We may further assume $\epsilon_2$ is
small so that $K_1+1-\frac{1}{\sqrt{2n\epsilon_2}}<0$.
Choose $p$ so that $2p(p-1)n\epsilon_2=1$, so
$2np^2\geq\frac{1}{\epsilon_2}$.  Therefore
 $K_1-p+2p(p-1)n\,\epsilon_2\leq K_1+1
 -\frac{1}{\sqrt{2n\epsilon_2}}<0 $.

Denote $\eta^{-2p}|A|^2 $ by $g$, then
$g$ satisfies
\[
\begin{split}
\frac{d}{dt}g
\leq \Delta g-\eta^{2p}
 \nabla(\eta^{-2p})\cdot
\nabla g+\eta^{2p}
(K_1+1-\frac{1}{\sqrt{2n\epsilon_2}})
g^2+K_2 g
\end{split}
\]

By the maximal principle and comparison theorem for parabolic
equations and notice that $0<\eta<1$ is bounded away from zero ,
$\max_{\Sigma_t} |A|^2\leq \frac{c_3 K_2}{\frac{1}
{\sqrt{2n\epsilon_2}}- (K_1+1)}$ if $t$ is large enough. Since
$\epsilon_2$ can be arbitrarily small, this implies
$\max_{\Sigma_t} |A|^2\rightarrow 0$ as $t\rightarrow \infty$.
Since the mean curvature flow is a gradient flow, the metrics are
analytic, by Simon's \cite{si} theorem, the flow converges to a
unique limit at infinity. That the limiting map is a constant map
follows from $*\Omega
=\frac{1}{\sqrt{\prod_{i=1}^n(1+\lambda_i^2)}}\rightarrow 1$ as
$t\rightarrow \infty$, thus $\lambda_i\rightarrow 0$ and
$|df_t|\rightarrow 0$.
\end{proof}


\begin{thebibliography}{99}
\bibitem{eh} K. Ecker and G. Huisken,
\textit{ Interior estimates for hypersurfaces
moving by mean curvature.}, Invent. Math. 105
(1991), no. 3, 547--569.

\bibitem{es} J. Eells and J. H. Sampson,
\textit{Harmonic mappings of Riemannian
manifolds}, Amer. J. Math. \textbf{86} (1964) 109--160.


\bibitem{hu2} G. Huisken, \textit{Asymptotic
behavior for singularities of the mean curvature flow}, J.
Differential Geom. \textbf{31} (1990), no. 1, 285--299.

\bibitem{il} T. Ilmanen, \textit{Singularities
of mean curvature flow of surfaces }, preprint , 1997. Available
at http://www.math.ethz.ch/~ilmanen/papers/pub.html

\bibitem{leon} S.J. Leon, \textit{Linear algebra with
applications}, 5th ed., Prentice Hall, Inc., 1998.


\bibitem{sch} R. Schoen, \textit{
The role of harmonic mappings in rigidity and
 deformation problems}, Complex geometry
 (Osaka, 1990), 179--200, Lecture Notes in
 Pure and Appl. Math., 143, Dekker, New York,
 1993.

\bibitem{si} L. Simon, \textit{Asymptotics for
 a class of nonlinear evolution equations,
 with applications to geometric problems.},
  Ann. of Math. (2) 118 (1983), no. 3, 525--571.

\bibitem{ke} M-T. Wang :\textit{Mean curvature flow of surfaces in Einstein
Four-Manifolds}, to appear in J. Differential Geom.

\bibitem{lg} M-T. Wang : \textit{Deforming area preserving diffeomorphism
of surfaces by mean curvature flow}, to appear in Math. Res. Lett.


\bibitem{w} B. White, \textit{A local
regularity theorem for classical mean curvature
flow}, preprint, 2000.



\end{thebibliography}
\end{document}